\documentclass{article}

\usepackage{authblk}

\bibliography{localbibliography}

\author[1]{Sebastien Bonduelle}
\author[2,3]{Franti\v{s}ek Kardo\v{s}}

\affil[1]{\small ENS de Rennes, France}
\affil[2]{\small Comenius University, Bratislava, Slovakia}
\affil[3]{\small Université de Bordeaux, CNRS,  LaBRI,  F-33400 Talence, France {\tt sebastien.bonduelle@ens-rennes.fr, frantisek.kardos@u-bordeaux.fr}}

\title{Subcubic planar graphs of girth 7 are class I}

\usepackage{geometry}
\usepackage{graphicx}
\usepackage{amsthm}
\usepackage{amsmath}
\usepackage{amssymb}
\usepackage{enumitem}
\usepackage{caption}
\usepackage{subcaption}
\usepackage{authblk}
\usepackage{tcolorbox}
\usepackage{algorithm2e}
\usepackage{hyperref}

\RestyleAlgo{ruled}
\SetKw{Not}{not}
\SetKw{Yield}{yield}
\SetKw{fYield}{yield from}
\SetKw{Break}{break}

\newtheorem{theorem}{Theorem}
\newtheorem{claim}{Claim}
\newtheorem{conjecture}{Conjecture}

\theoremstyle{definition}
\newtheorem{definition}{Definition}[section]

\theoremstyle{remark}

\newcommand{\degree}{\mathrm{deg}}

\begin{document}
\maketitle

\begin{abstract}
We prove that planar graphs of maximum degree 3 and of girth at least 7 are 3-edge-colorable, extending the previous result for girth at least 8 by Kronk, Radlowski, and Franen from 1974. 
\end{abstract}

\section{Introduction}

One of the most famous results in graph theory, and the most prominent result concerning edge coloring of graphs, is Vizing's theorem, stating that every (simple) graph of maximum degree $\Delta$ admits an edge-coloring using at most $\Delta+1$ colors, thus reducing the possible values of the chromatic index of a graph to only $\Delta$ and $\Delta+1$ \cite{Viz}. 
Graphs admitting a $\Delta$-edge-coloring are called class I, those who need $\Delta+1$ colors class II.

Vizing also proved that, among other classes, planar graphs of maximum degree at least 8 are class I \cite{Viz2}.
Kronk, Radlowski, and Franen \cite{KRF} proved that planar graphs of girth at least $g$ and of maximum degree at least $d$ are class I for $(g,d)=(3,8), (4,5), (5,4), (8,3)$, generalizing Vizing's result. See also the paper by Li and Luo \cite{LiLuo} where these results are generalized to other surfaces and reproved using modern language.

On the other hand, examples of planar graphs of class II for $(g,d)=(3,5), (3,4), (3,3), (4,3)$, and $(5,3)$ have been known for a long time: it suffices to consider a class I regular graph (e.g., a graph of a platonic solid) and subdivide one edge.

This left only five cases open: $(g,d)=(3,7), (3,6), (4,4), (6,3), (7,3)$.

Sanders and Zhao \cite{SZh}, and independently Zhang \cite{Zhang}, proved that planar graphs of maximum degree 7 are class I, closing therefore the first open case. 
There has been a lot of attention paid to the second open case, see the survey \cite{survey} for further results.

In this paper, we close the fifth open case in the affirmative. We prove
\begin{theorem}
Planar subcubic graphs of girth at least 7 are class I.
\label{th:main}
\end{theorem}

In the rest of the paper, whenever we speak about a coloring, we mean a proper 3-edge-coloring using the colors from $\{1, 2, 3\}$.

All the graphs considered here are simple (i.e., without loops and multiple edges). If we speak about a pseudograph, then multiple edges and loops are allowed; however, pseudographs appearing in the paper do not have multiple edges.

We will distinguish between a \emph{planar graph} (a graph such that there exists a planar embedding without crossing edges) and a \emph{plane graph} (an actual crossing-free embedding of a planar graph in the plane). 

\section{Proof of Theorem \ref{th:main}}

Let $G_0$ be a counterexample with the smallest number of vertices, and subject to that, the smallest number of edges. Clearly, we may assume that $G_0$ is 2-connected (and, in particular, it has no vertices of degree at most one), and that it has girth exactly 7. By minimality, every graph obtained from $G_0$ by removing a nonempty set of edges and/or vertices, is 3-edge-colorable. We fix a plane embedding of $G_0$.

\subsection{Reducible patterns}

Let $G$ be a 2-connected plane graph. Let $F=(X,Y)$, with $F\subset E(G)$, $X\cup Y = V(G)$, and $X\cap Y = \emptyset$, be an inclusion-wise minimal edge-cut in $G$ such that the outer face of $G$ is incident with at least one vertex from $Y$. By minimality, $G[X]$ and $G[Y]$ are connected. 
Let $P_X$ be the graph obtained from $G[X]$ by adding a pending half-edge into the outer face for each cut edge from $F$. Observe that $\degree_{P_X}(v)=\degree_G(v)$ for all $v\in X$.

\begin{definition}[Pattern]
A \emph{pattern} is a connected plane graph $H$ in which vertices incident with the outer face may have pending half-edges towards the outer face. We say that a pattern $H$ is \emph{contained} in a 2-connected plane graph $G$, if $H$ is isomorphic to $P_X$ for some inclusion-wise minimal edge-cut $(X,Y)$ of $G$. If this is the case, we call $P_X$ a \emph{realization} of $H$ in $G$.
\end{definition}

Here we will only consider subcubic patterns, with an additional property that every vertex has at most one incident half-edge.

\begin{definition}[Frontier coloring]
    Let $H$ be a pattern. The \emph{frontier} of $H$, denoted by $\partial H$, is the set of its half-edges, called \emph{frontier edges}. 
    A \emph{frontier coloring} of $H$ is an assignment of colors from $\{1, 2, 3\}$ to the elements of $\partial H$. We denote by $\Gamma(H)= \{1, 2, 3\}^{\partial H}$ the set of frontier colorings of $H$.
\end{definition}

    If $\varphi:E(H)\to \{1,2,3\}$ is a coloring of a pattern $H$, then $\varphi|_{\partial H}$ is a frontier coloring. On the other hand, not necessarily every frontier coloring of $H$ can be extended into a coloring of $H$. We denote by $\Gamma_0(H)$ those that are extendable.

\begin{definition}[Boundary switch]
Let $\gamma$ be a frontier coloring of a pattern $H$, let $\{i,j\}$ be a pair from \{1, 2, 3\}, and let $e$ be a frontier edge colored $i$ or $j$ by $\gamma$. We denote $\gamma_{ij}(e)$ the frontier coloring obtained from $\gamma$ by switching the color of $e$ from $i$ to $j$ or vice versa. Similarly, for a pair of edges $e,e'$ (with $e\ne e'$) colored $i$ or $j$ (not necessarily both with the same color), we denote $\gamma_{ij}(e,e')$ the frontier coloring obtained from $\gamma$ by switching the colors of $e$ and $e'$ from $i$ to $j$ or vice versa.
\end{definition}

\begin{definition}[Auxiliary graph with respect to a set of frontier colorings and a color pair]
    Let $H$ be a pattern and $\gamma$ a frontier coloring of $H$. Let $\Gamma$ be a subset of $\Gamma(H)$ and let $\{i,j\}$ be a pair from \{1, 2, 3\}. The \emph{auxiliary graph of $\gamma$ with respect to $\Gamma$ and $\{i,j\}$} is the pseudograph embedded in the plane (with possibly some edges crossing) denoted by $\mathcal{A}_{H, \Gamma, \gamma}(i, j)$ where:
    
    \begin{itemize}
        \item The vertices of $\mathcal{A}_{H, \Gamma, \gamma}(i, j)$ represent the frontier edges of $H$ that are colored $i$ or $j$ by $\gamma$. They are disposed on a circle $K$ in the plane respecting the cyclic order of the corresponding frontier edges, thus forming vertices of a convex polygon.
        \item There is a straight-line edge between vertices $e$ and $e'$ (a loop at $e$ outside $K$) if $\gamma_{ij}(e,e') \notin \Gamma$ ($\gamma_{ij}(e) \notin \Gamma$, respectively).
    \end{itemize}
\end{definition}

\begin{definition}[Non-crossing perfect quasi-matching]
 A \emph{quasi-matching} of a pseudograph $\mathcal{A}$ is a set of pairwise disjoint edges or loops of $\mathcal{A}$. A quasi-matching $M$ is \emph{perfect} if it is spanning, i.e. every vertex of $\mathcal{A}$ is incident to at least (and thus, exactly) one element of $M$. If a plane embedding of $\mathcal{A}$ is given, a quasi-matching $M$ is \emph{non-crossing} if no two edges of $M$ do cross.
\end{definition}

For a realization of a pattern $H$ contained in a plane graph $G$, a non-crossing perfect quasi-matching in $\mathcal{A}_{H, \Gamma, \gamma}(i, j)$ would correspond to a particular position of $(i,j)$-Kempe chains in a coloring of $(G\setminus H)\cup \partial H$ such that no single $(i,j)$-Kempe switch leads to a frontier coloring in $\Gamma$.

\begin{definition}[Reduction to a set of frontier colorings]
    Let $H$ be a pattern, $\gamma \in \Gamma(H)$ and $\Gamma \subseteq \Gamma(H)$. The coloring $\gamma$ is \emph{reducible to $\Gamma$} if there exists a pair $i \neq j$ from \{1, 2, 3\} such that $\mathcal{A}_{H, \Gamma, \gamma}(i, j)$ has no non-crossing perfect quasi-matching.
\end{definition}

A coloring $\gamma$ being reducible to $\Gamma$ means that for a certain color pair $(i, j)$, whichever way the $(i, j)$-Kempe chains starting from the frontier edges are disposed in, it will be possible to swap one of them to reach a frontier coloring in $\Gamma$.

\begin{definition}[Rank of a frontier coloring]
    We define inductively the \emph{rank} of a frontier coloring as follows:  
    \begin{itemize}[label=$\bullet$]
        \item $\gamma$ is of rank 0 if $\gamma \in \Gamma_0(H)$;
        \item $\gamma$ is of rank $k + 1$ if it is not in $\Pi_k(H)$ and $\gamma$ is reducible to $\Pi_k(H)$,
    \end{itemize}
    where $\Pi_k(H)$ denotes the set of frontier colorings of rank up to $k$. 
\end{definition}
    We denote by
    $\Pi(H) = \bigcup_{k=0}^\infty \Pi_k(H)
    $ the set of all frontier colorings that have a rank.

By finiteness of the set $\Gamma(H)$ and by monotony of $(\Pi_k(H))_{k=0}^\infty$, there exists $k_0$ such that $\Pi_k(H)=\Pi_{k_0}(H)$ for all $k\ge k_0$. Additionally, the minimal $k$ such that $\Pi_{k + 1}(H) = \Pi_{k}(H)$ is such a $k_0$. Indeed, if
$\Pi_{k + 1}(H) = \Pi_k(H)$, the frontier colorings of rank $k + 2$ are exactly those of rank $k + 1$, of which there are none. So $\Pi_{k + 2}(H) = \Pi_{k + 1}(H)$, and we proceed by recursion to conclude. This means that once we find a rank empty of colorings, we can stop and know that all the remaining colorings will not have a rank.

\begin{definition}[Reducibility of a pattern]
    A pattern $H$ is \emph{reducible} if $\Pi(H)=\Gamma(H)$.
\end{definition}

In other words, a pattern $H$ is reducible if, for every graph $G$ containing $H$ and for every coloring $\varphi$ of $G\setminus H \cup \partial H$,  the frontier coloring induced by $\varphi_{|\partial H}$ either directly extends to a coloring of $H$ (and thus a coloring of $G$), or there exists a finite sequence of Kempe switches that transforms $\varphi$ into an extendable one. Clearly, $G_0$ contains no reducible pattern.

This definition generalizes the notion of $D$-reducibility, used in the proof of the Four Color Theorem \cite{4CT}, see also \cite{Stein}. The main difference is that here we are dealing with subcubic and not cubic graphs, and so we cannot make use of the parity lemma. Furthermore, we cannot simply look at non-crossing perfect matchings in the auxiliary graph, we have to take loops into consideration as well.

\begin{claim}
The following patterns are reducible:
\begin{itemize}
    \item $P_{2^2}$, a path of two vertices of degree 2;
    \item $P_{232}$, a path of three vertices of degrees 2,3,2, respectively;
    \item $P_{3_2^3}$, a path of three vertices of degree 3, each having a neighbor of degree 2, all on the same side of the path;
    \item $P_{233_2^2}$, a path of three vertices of degree 3, each having a neighbor of degree 2, the first on one side of the path, the other two on the other side;
    \item $P_{3_2^233_2}$, a path of four vertices of degree 3, of which the first two and the last one have a neighbour of degree 2, all on the same side of the path;
    \item $P_{7}$, a pattern composed of one inner face of size 7, with exactly one vertex of degree 2; moreover each end-vertex of the opposite edge has a neighbor of degree 2.
\end{itemize}
\end{claim}

\begin{figure}[h!]
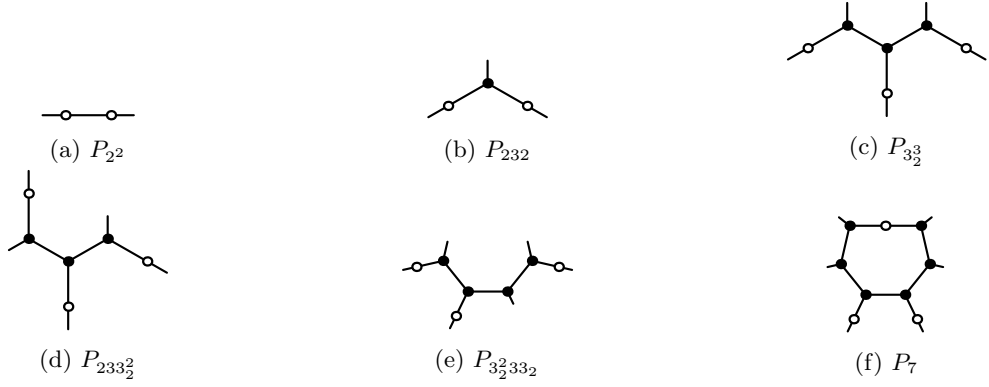

    \begin{subfigure}[b]{0.3\textwidth}
        \centering
        \includegraphics[scale=1]{./Figures/patternA.mps}
        \caption{$P_{2^2}$}
    \end{subfigure}
    \hfill
    \begin{subfigure}[b]{0.3\textwidth}
        \centering
        \includegraphics[scale=1]{./Figures/patternB.mps}
        \caption{$P_{232}$}
    \end{subfigure}
    \hfill
    \begin{subfigure}[b]{0.3\textwidth}
        \centering
        \includegraphics[scale=1]{./Figures/patternC.mps}
        \caption{$P_{3_2^3}$}
    \end{subfigure}
    \hfill
    \begin{subfigure}[b]{0.3\textwidth}
        \centering
        \includegraphics[scale=1]{./Figures/patternD.mps}
        \caption{$P_{233_2^2}$}
    \end{subfigure}
    \hfill
    \begin{subfigure}[b]{0.3\textwidth}
        \centering
        \includegraphics[scale=1]{./Figures/patternE.mps}
        \caption{$P_{3_2^233_2}$}
    \end{subfigure}
    \hfill
    \begin{subfigure}[b]{0.3\textwidth}
        \centering
        \includegraphics[scale=1]{./Figures/pattern7.mps}
        \caption{$P_7$}
    \end{subfigure}
    \caption{Reducible patterns used in the proof of Theorem \ref{th:main}. Vertices of degree 3 are represented by dots, those of degree 2 by empty circles.}
\label{fig:patterns}
\end{figure}

See Figure \ref{fig:patterns} for an illustration of these  patterns.

\begin{proof} We have implemented a program to check the reducibility of the patterns. Instead of considering all the frontier colorings, we only considered equivalence classes modulo color permutations and symmetries of the patterns. 
A hand-made case-analysis proof of the reducibility of each of the patterns, as well as the \href{https://github.com/asimov-io/pattern-reducibility-checker}{source codes}, are available on arxiv \cite{arxiv}.

The program takes as input a pattern $H$, represented by the adjacency list of its line graph $L(H)$, together with the list of its frontier edges in a cyclic order, and its symmetry group.

As a preprocessing stage, first, a map $R()$, mapping each frontier coloring $\gamma$ to its representative $R(\gamma)$ (the lexicographically minimal frontier coloring equivalent to $\gamma$ modulo color permutations and symmetries), is defined and computed. Then, the set $\bar{\Gamma}$ of equivalence class representatives is calculated.

The rank $r$ of each $R(\gamma)$ is initialized to $+\infty$. Using an exhaustive backtracking subroutine, $\Gamma_0(H)$ is computed -- the rank of $R(\gamma)$ is reset to 0 if (and only if) $\gamma$ is extendable into a 3-edge-coloration of $H$.

Then, in the main loop, the ranks of all the remaining reducible frontier colorings are determined. See Algorithm \ref{alg:redCheck} for further details.

\begin{algorithm*}
    \caption{Reducibility Checker}
    \label{alg:redCheck}

    Compute $R(\gamma)$ for each $\gamma \in \Gamma(H)$
    
    Compute $\Gamma_0(H)$: set $r[R(\gamma)]=0$ for each $\gamma \in \Gamma_0(H)$, set $r[R(\gamma)]=+\infty$ otherwise
    
    $i = 1$\\
    
    $found\_changed = True$\\
    
    \While{$found\_changed$}{
        $found\_changed = False$\\
    
        $found\_non\_reducible = False$\\
    
        \For{$\gamma \in \bar{\Gamma}$}{
            \If{$r[\gamma]=+\infty$}{

                \For{$k \in \{1, 2, 3\}$}{
                    Build the auxiliary graph of $\gamma$ with respect to $\{\gamma' \in \Gamma(H) \mid r(R(\gamma'))<i\}$ and $\{1, 2, 3\} \setminus \{k\}$\\
                    then test if it admits a non-crossing perfect quasi-matching
                }
                \eIf{\emph{each auxiliary graph admits a matching}}{
                    $found\_non\_reducible = True$\\
                }{
                    $r[\gamma]=i$;\\
                    $found\_changed = True$\\
    
                }
            }
        }
    
        $i += 1$\\
    }
    \Return{$found\_non\_reducible$, $r[\,]$}

\end{algorithm*}

In Table \ref{tab:red} we provide an overview of the numbers of equivalence classes of colorings of rank $k$ for relevant values of $k$, for our six reducible patterns.

\begin{table}[ht]
    \centering
    \begin{tabular}{|c||c||c|c|c|c|c|c|}\hline
         Pattern & total & rank 0 & rank 1 & rank 2 & rank 3 & rank 4 & rank 5 \\\hline\hline
         $P_{2^2}$ & 2 & 2 & -- & -- & -- & -- & --  \\\hline
         $P_{232}$ & 4 & 3 & 1 & -- & -- & -- & --\\\hline
         $P_{3_2^3}$ & 25 & 14 & 5 & 4 & 2 & -- & --\\\hline
         $P_{233_2^2}$ & 41 & 26 & 9 & 5 & 1 & -- & --\\\hline
         $P_{3_2^233_2}$ & 122 & 56 & 23 & 13 & 14 & 15 & 1 \\\hline
         $P_{7}$ & 70 & 38 & 13 & 12 & 5 & 2 & -- \\\hline
    \end{tabular}
    \caption{For each of the reducible patterns, we provide the number of equivalence classes of frontier colorings (with respect to permutations of colors and to symmetries of the pattern), along with the numbers of colorings of rank $k$.}
    \label{tab:red}
\end{table}

See Figure \ref{fig:232} for the illustration of all possible (equivalence classes of) colorings and their rank for the pattern $P_{232}$.

\begin{figure}[h!]
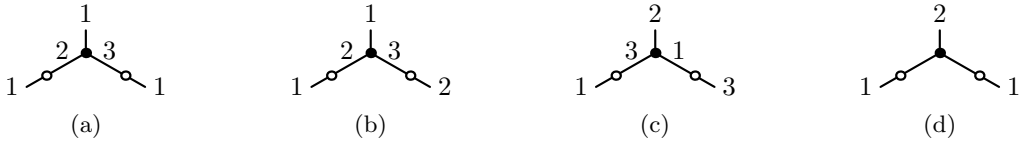

    \centerline{
    \begin{subfigure}[b]{0.24\textwidth}
        \centering
        \includegraphics{Figures/pattern232a.mps}
        \caption{}
    \end{subfigure}
    \hfil
    \begin{subfigure}[b]{0.24\textwidth}
        \centering
        \includegraphics{Figures/pattern232b.mps}
        \caption{}
    \end{subfigure}
    \hfil
    \begin{subfigure}[b]{0.24\textwidth}
        \centering
        \includegraphics{Figures/pattern232d.mps}
        \caption{}
    \end{subfigure}
    \hfil
    \begin{subfigure}[b]{0.24\textwidth}
        \centering
        \includegraphics{Figures/pattern232c.mps}
        \caption{}
    \end{subfigure}}
    \caption{For the pattern $P_{232}$, there are only four equivalence classes of frontier colorings. Three of them are of rank 0, the last one reduces to the previous three: in (d), for the pair of colors (1, 2), if an edge colored 1 is switched, the coloring class (b) is obtained. If the middle edge is switched, the coloring class (a) is obtained. Therefore, the auxiliary graph contains no loops; and since it has an odd number of vertices, it cannot admit a perfect quasi-matching.}
    \label{fig:232}
\end{figure}

\end{proof}

\subsection{Discharging}

We now know that the counterexample $G_0$ does not contain any of those patterns because they are reducible. In particular, the vertices of degree 2 are at distance at least 3 from each other.

From Euler’s formula, we have $6V+6F-6E = 12 > 0$, and so \[\sum_{v \in V} (6-2 \degree(v)) + \sum_{f \in F} (6-\degree(f)) > 0.\]

We assign a charge of $6-2\degree(v)$ to a vertex $v$ and a charge of $6-\degree(f)$ to a face $f$. Clearly, the sum of all charges is positive.
Note that vertices of degree 2 are the only elements of the graph that contain positive charge; all the faces are negative.

We apply the two following discharging rules, one after another, in order to observe the annihilation of all positive charges, and thus obtain a contradiction.
\begin{enumerate}
    \item [(R1)] Every vertex of degree 2 sends $1$ of charge to each incident face.
\end{enumerate}

\begin{figure}[h!]
    \centering
    \includegraphics{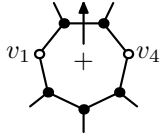}
    \caption{(R2) applied to a critical face.}
    \label{fig:R2}
\end{figure}

Observe that after applying (R1) all the positive charge has been moved to faces. Since the distance between any two vertices of degree 2 in $G_0$ is at least 3 (because $P_{2^2}$ and $P_{232}$ are reducible), a face of size $d$ receives at most $\lfloor \frac{d}{3}\rfloor$ units of charge. If $d\ge 9$, then its charge is at most $6-d +\frac{d}{3} = \frac23 (9-d)\le 0$. If $d=8$, then its charge is at most $6-d+\lfloor \frac{d}{3}\rfloor = 6-8+2 = 0$. 
Therefore, the only element that can have positive charge after (R1) is a face of size 7, incident with two vertices of degree 2 -- it will have a charge of $+1$. We will call such a face a \emph{critical} face. We now apply the following rule:
\begin{enumerate}
    \item [(R2)] Let $f$ be a critical face. Let $v_1v_2v_3v_4v_5v_6v_7$ be the facial cycle of $f$, where $v_1$ and $v_4$ are the two vertices of degree 2 incident to $f$. Let $f'$ be the face adjacent to $f$ via the edge $v_2v_3$. Then $f$ sends its charge (of $1$) to the face $f'$ (see Figure \ref{fig:R2}).
\end{enumerate}

It is clear that after (R2) all the vertices remain free of charge. It remains to prove that no new face with positive charge is created. 

Let $f$ be a face of size $d\ge 7$. Its initial charge is $6-d$. It can receive charge from incident vertices of degree 2, and from adjacent critical faces.

\begin{figure}[h!]
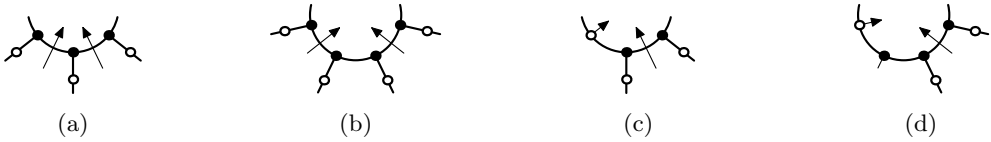

    \begin{subfigure}[b]{0.24\textwidth}
        \centering
        \includegraphics{./Figures/pattern_critical2_dist_1.mps}
        \caption{}
        \label{fig:crit^2_dist_1}
    \end{subfigure}
    \begin{subfigure}[b]{0.24\textwidth}
        \centering
        \includegraphics{./Figures/pattern_critical2_dist_2.mps}
        \caption{}
        \label{fig:crit^2_dist_2}
    \end{subfigure}
    \begin{subfigure}[b]{0.24\textwidth}
        \centering
        \includegraphics{./Figures/pattern_2-critical_dist_1.mps}
        \caption{}
        \label{fig:2-crit_dist_1}
    \end{subfigure}
    \begin{subfigure}[b]{0.24\textwidth}
        \centering
        \includegraphics{./Figures/pattern_2-critical_dist_2.mps}
        \caption{}
        \label{fig:2-crit_dist_2}
    \end{subfigure}
    \caption{If two critical faces are at distance 1 or 2 around a face, then the graph contains a reducible pattern $P_{3_2^3}$ (\ref{fig:crit^2_dist_1}-\ref{fig:crit^2_dist_2}). If a vertex of degree 2 and a critical face are at distance 1.5 (2.5) around a face, then the graph contains a reducible pattern $P_{232}$ (\ref{fig:2-crit_dist_1}) ($P_{233_2^2}$ (\ref{fig:2-crit_dist_2}), respectively).}
    
    \label{fig:dist_3}
\end{figure}

Two elements sending charge to $f$ must be at distance at least 3 along the boundary of $f$, otherwise a reducible pattern $P_{232}$, $P_{3_2^3}$ or  $P_{233_2^2}$ can be found in $G_0$ (see Figure \ref{fig:dist_3}), a contradiction.

Therefore, $f$ receives at most $\lfloor \frac{d}{3}\rfloor$ units of charge by (R1) and (R2) in total. Again, a face can only become positive if $d=7$ and it receives 2 units of charge. If it is due to two vertices of degree 2, then $f$ is a critical face and it evacuates its charge by (R2). If $f$ receives charge from two critical faces, or from a critical face and a vertex of degree 2, then $G_0$ contains a reducible pattern $P_{3_2^233_2}$ or $P_{7}$, respectively, a contradiction.

This concludes the proof of Theorem \ref{th:main}.

\section{Concluding remarks}

Note that the edge-coloring problem is known to be NP-complete \cite{Hoyler}, and it remains NP-hard even for triangle-free subcubic graphs \cite{Koreas}.

We recall a conjecture of unclear origin, often attributed to Vizing:
\begin{conjecture}
Subcubic planar graphs of girth 6 are class I.
\end{conjecture}

Considering the reducibility as defined in this paper does not seem to be sufficient to find a set of unavoidable reducible patterns in order to prove this conjecture. In particular, we haven't found any reducible pattern containing only one vertex of degree 2, which could be necessary since faces of size 6 incident to one vertex of degree 2 could become positive without having a negative element at arbitrary large neighborhood to discharge them.

The Four Color Theorem is equivalent to the statement that cubic planar graphs of girth 5 are class I. 
It would be interesting to see whether the conjecture above implies this statement.
\bigskip

When computing the reducibility of the patterns, we have encountered the following problem:

{\sl Does a given pseudo-graph with a fixed plane embedding (allowing edge crossings) admit a non-crossing perfect quasi-matching?}

For the goal of this work, we have polynomially reduced this problem to CNF-SAT (see \href{https://github.com/asimov-io/pattern-reducibility-checker/blob/main/aux_graph.py}{aux\_graph.py}) and used a SAT-solver in our program --- for us it was sufficient because the instances considered here were very small (pseudo-graphs with at most 6 vertices).

It could be interesting to determine if the non-crossing perfect quasi-matching problem is NP-hard, or if there exists a better algorithm than to reduce it to CNF-SAT.

\subsubsection*{Acknowledgement}
This work was supported in part by APVV-19-0308 and by VEGA 1/0743/21.

\vfill\eject 
\section*{Annex A: Proof of reducibility of the patterns without the help of a computer}

The reducibility of the patterns $P_{2^2}$ and $P_{232}$ was already known to Vizing.

Before considering the remaining patterns one by one, we introduce some common properties of patterns containing a facial path of four vertices with degrees $2,3,3,2$.

\subsection*{Common observations}

Let $P$ be a pattern containing a facial path $u'uvv'$ of vertices of degrees $2,3,3,2$, respectively. Let $e=uv$, $e_{u'}=uu'$, $e_{v'}=vv'$; let $f_{u'}$ ($f_{v'}$) be the other edge incident with $u'$ ($v'$, respectively), let $e_u$ ($e_v$) be the third edge incident with $u$ ($v$, respectively) not on the path.

Let $G'=G\setminus e$; let $\varphi$ be a coloring of $G'$. We have $\{1,2,3\}=\{\varphi(e_{u}),\varphi(e_{u'}), \varphi(e_v),\varphi(e_{v'})\}$, since otherwise $\varphi$ extends to a coloring of $G$ easily.

Without loss of generality we may assume that the edges incident with $u$ are colored with colors 1 and 2, and that the edges incident with $v$ are colored with colors 1 and 3.

We claim that we can reduce all the cases to the case where $\varphi(e_{u'})=\varphi(e_{v'})=1$.

\begin{figure}[h!]
    \centerline{
    \includegraphics{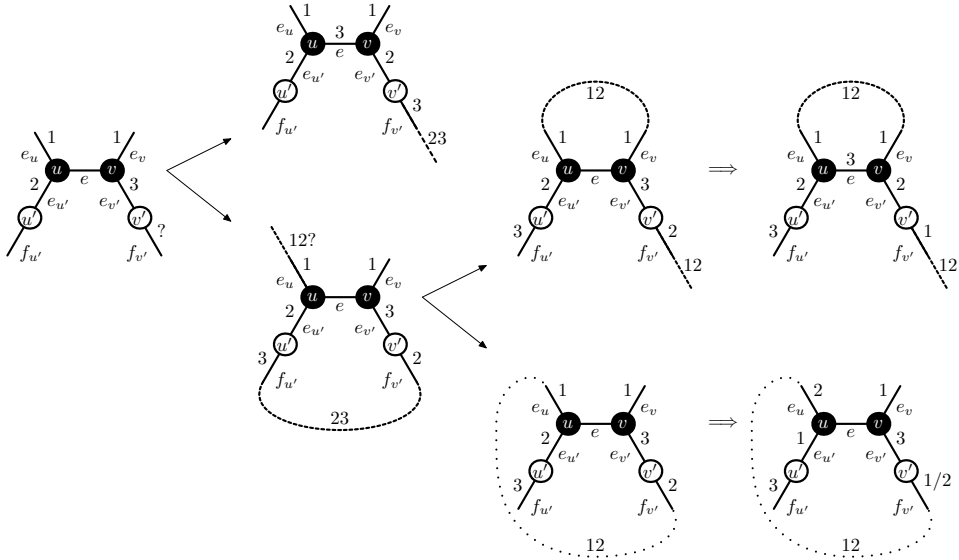}
    }
    \caption{The case when $\varphi(e_u)=\varphi(e_v)=1$ can be reduced to the case when $\varphi(e_{u'})=\varphi(e_v)=1$. Dashed lines represent confirmed Kempe chains, dotted lines represent potential Kempe chains.}
    \label{fig:1213}
\end{figure}

Suppose first that neither $e_{u'}$ nor $e_{v'}$ is colored 1. Then there is a (2,3)-Kempe chain from $u$ to $v$ starting with $e_{u'}$ and ending with $e_{v'}$ (otherwise we could recolor one of them to the color of the other and complete the coloring). In particular, we have $\varphi(f_{u'})=3$ and $\varphi(f_{v'})=2$. Therefore, there is a (1,2)-Kempe chain starting with $e_{u'}$ at $u'$. If this chain ends with $e_v$ at $v$, then $f_{v'}$ does not belong to this chain, and we may switch it to 1, recolor $e_{v'}$ to 2 and complete the coloring by setting $\varphi(e)=3$. Hence, the (1,2)-Kempe chain starting with $e_{u'}$ at $u'$ either does not come back to the pattern at all, or it ends with $f_{v'}$ at $v'$. After switching this Kempe chain, we have $\varphi(e_{u'})=1$ and $\varphi(e_{v'})=3$, the color of $f_{v'}$ may have been switched as well or not. See Figure \ref{fig:1213} for illustrations.

\begin{figure}[h!]
    \centerline{
    \includegraphics{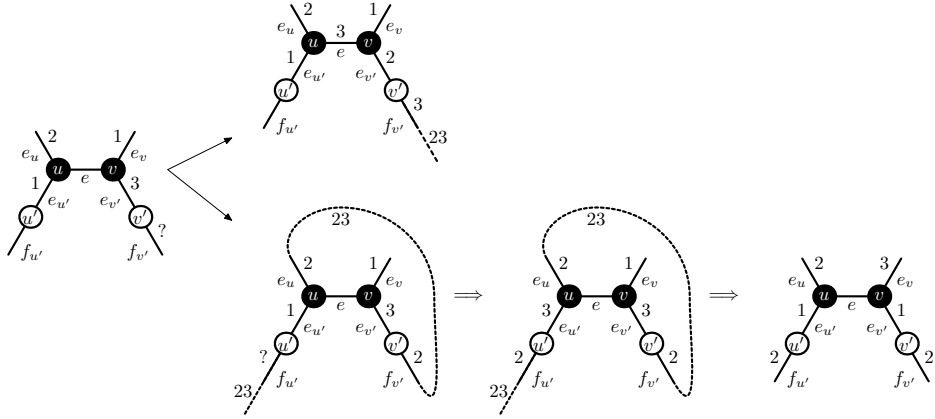}
    }
    \caption{The case when $\varphi(e_u')=\varphi(e_v)=1$ can be reduced to the case when $\varphi(e_{u'})=\varphi(e_{v'})=1$. }
    \label{fig:2113}
\end{figure}

Suppose now that one of $e_{u'}$ and $e_{v'}$ is colored 1. Without loss of generality we may assume that $\varphi(e_{u'})=1$ and $\varphi(e_{v'})=3$, and so $\varphi(e_u)=2$ and $\varphi(e_v)=1$. Then there is a (2,3)-Kempe chain from $u$ to $v$ starting with $e_u$ and ending with $e_{v'}$; in particular, $\varphi(f_{v'})=2$.
The (2,3)-Kempe chain starting with $f_{u'}$ at $u'$ is disjoint from the previous one, and so we may switch it freely to get $\varphi(f_{u'})=2$.
After the switch we can recolor $e_{u'}$ with 3, and so (after a color permutation) we obtain a coloring with the desired property. See Figure \ref{fig:2113} for illustration.

From now on we may assume that $\varphi(e_{u'})=\varphi(e_{v'})=1$ and that $\varphi(e_u)=2$ and $\varphi(e_v)=3$. Moreover, there is a (2,3)-Kempe chain from $u$ to $v$ starting with $e_{u}$ and ending with $e_{v}$. Two other (2,3)-Kempe chains start with $f_{u'}$ at $u'$ and with $f_{v'}$ at $v'$. If these are two distinct Kempe chains, then we can switch them independently from each other to obtain $\varphi(f_{u'})=2$ and $\varphi(f_{v'})=3$, and complete the coloring by setting $\varphi(e)=1$, $\varphi(e_{u'})=3$ and $\varphi(e_{v'})=2$. Hence, we may assume that the (2,3)-Kempe chain that starts with $f_{u'}$ at $u'$ ends with $f_{v'}$ at $v'$, moreover, the colors of $f_{u'}$ and $f_{v'}$ are the same.

\subsection*{The patterns $P_{3_2^3}$ and $P_{233_2^2}$}

We will deal with both patterns at the same time. We introduce a notation that applies to both of them.

Let $u'uvv'$ be a facial path of four vertices of degrees $2,3,3,2$, with $e$, $e_u$, $e_v$, $e_{u'}$, $e_{v'}$, $f_{u'}$, and $f_{v'}$ defined as before. Let $e_v=vw$ and suppose $w$ has a neighbor of degree 2, say $w'$. 
Let $e_{w'}=ww'$ and let $e_w$ ($f_{w'}$) be the frontier edge incident with $w$ ($w'$, respectively). 

According to the common observations, we may assume that $\varphi(e_{u'})=\varphi(e_{v'})=1$, $\varphi(e_u)=2$ and $\varphi(e_v)=3$, $\varphi(f_{u'})=\varphi(f_{v'})=2$ and there are (2,3)-Kempe chains from $u$ to $v$ and from $u'$ to $v'$. 

\begin{figure}[h!]
    \centerline{
    \includegraphics{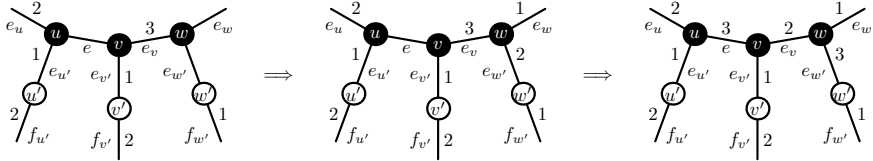}
    }
    \caption{A coloring with $\varphi(f_{w'})=1$ can be easily extended to the whole graph.}
    \label{fig:case1}
\end{figure}

\begin{figure}[h!]
    \centerline{
    \includegraphics{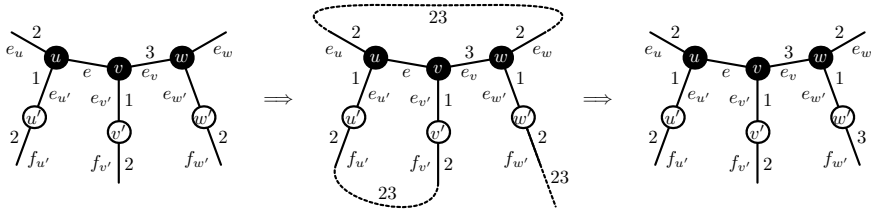}
    }
    \caption{A coloring with $\varphi(f_{w'})=2$ can be reduced to the case $\varphi(f_{w'})=3$.}
    \label{fig:case2}
\end{figure}

We claim that we may suppose that $\varphi(f_{w'})=3$: It is easy to see that $\varphi(f_{w'})\ne 1$, as otherwise we would have $\varphi(e_{w'})=2$ and the (2,3)-Kempe chain from $v$ would not continue, see Figure \ref{fig:case1}. If we had $\varphi(f_{w'})=2$, then a new (2,3)-Kempe chain (distinct from the two already established since $\varphi(e_{w'}) = 1$ in this case) would start with $f_{w'}$ at $w'$, and we could switch it, see Figure \ref{fig:case2} for illustration.


The (1,2)-Kempe chain leaving the pattern starting with $e_u$ at $u$ has to end with $f_{v'}$ and $e_{v'}$ at $v$, otherwise we could recolor $e_{u'}$ to 3, color $e$ with 2 and switch the (1,2)-Kempe chain starting with $e_u$ to obtain a coloring of $G$.

Therefore, the (1,2)-Kempe chain starting with $e_{w'}$ at $w'$ is distinct from the one connecting $e_u$ with $f_{v'}$. By planarity, this new chain cannot get to $f_{u'}$, and so we can switch it freely to obtain a coloring where $\varphi(e_w)=2$.

We are in the situation where four out of the five frontier edges are colored 2, and $f_{w'}$ is colored 3. The position of the (2,3)-Kempe chains may have changed, since we may have performed a (1,2)-switch.

\begin{figure}[h!]
    \centerline{
    \includegraphics{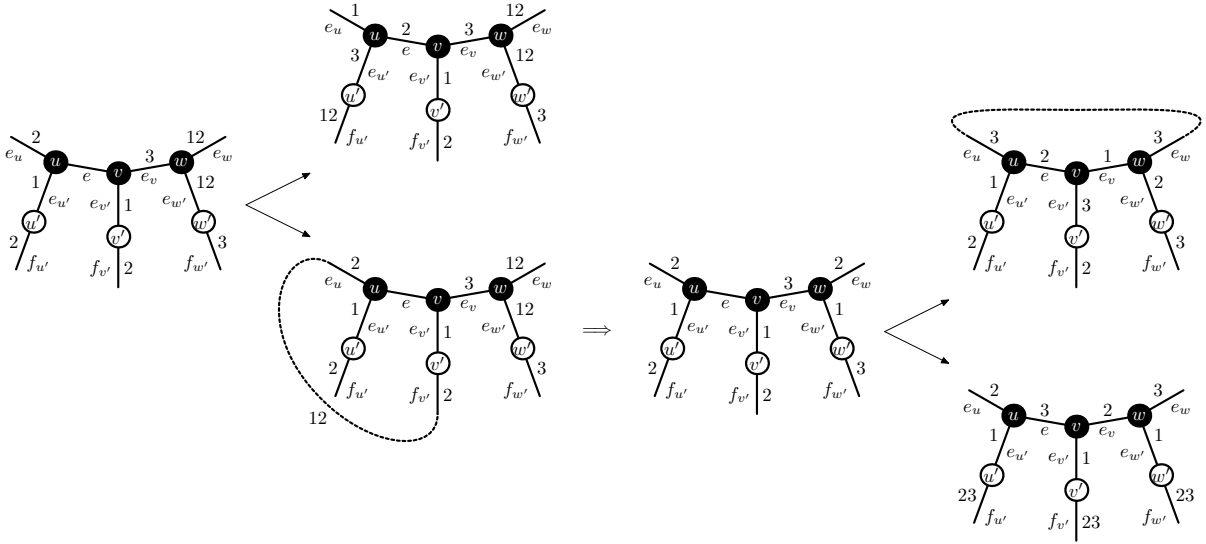}
    }
        \caption{A coloring with $\varphi(f_{w'})=3$ is reducible.}
    \label{fig:case3}
\end{figure}

We switch the (2,3)-Kempe chain starting with $e_v$ at $v$, recoloring $e_w$ to 3 on the way. If $e_u$ gets recolored to 3, then we recolor $e_{v'}$ to 3, $e_v$ to 1, $e_{w'}$ to 2, and we color $e$ with 2. If $e_u$ is not affected, then we simply color $e$ with 3, no matter whether any among $f_{u'}$, $f_{v'}$, or $f_{w'}$ gets switched from 2 to 3 or vice versa. See Figure \ref{fig:case3} for illustration.

\subsection*{The pattern $P_{3_2^2 33_2}$}

Let $u_1u_2u_3u_4$ be four consecutive vertices of degree 3 on a facial path, let $v_1$, $v_2$, and $v_4$ be vertices of degree 2 adjacent to $u_1$, $u_2$, and $u_4$, respectively. Let $w_i$ be the other neighbor (outside the pattern) of $v_i$, for $i=1,2,4$.
Let $u_0$ be the third neighbor of $u_1$, $v_3$ be the third neighbour of $u_3$ and $u_5$ be the third neighbor of $u_4$, all the three outside the pattern.

Let $\varphi$ be a coloring of $G\setminus u_1u_2$. From the common observations we may assume that $\varphi(u_1v_1)=\varphi(u_2v_2)=1$, $\varphi(u_0u_1)=3$, $\varphi(u_2u_3)=2$, and there is a (2,3)-Kempe chain starting with $u_1u_0$ at $u_1$ and ending with $u_3u_2$ at $u_2$.
Moreover, there is another (2,3)-Kempe chain starting with $v_1w_1$ at $v_1$ and ending with $w_2v_2$ at $v_2$, disjoint from the previous one. After an eventual switch, we may assume that $\varphi(v_1w_1)=\varphi(v_2w_2)=2$.

Furthermore, there is a (1,2)-Kempe chain starting with $v_1w_1$ at $v_1$ and ending with $u_3u_2$ at $u_2$.

The (1,2)- and (2,3)-Kempe chains starting with $u_2u_3$ at $u_2$ have to leave the pattern by either $u_3v_3$, $u_4u_5$, or $v_4w_4$.
At least one of them does not use the edge $u_3v_3$.

Suppose first that the (2,3)-Kempe chain starting with $u_2u_3$ at $u_2$ leaves the pattern by $u_4u_5$. Then there is a third (2,3)-Kempe chain starting with $v_4w_4$ at $v_4$ (because in this case $\varphi(u_4v_4) = 1$); after an eventual switch of that chain we obtain a coloring with $\varphi(v_4w_4)=2$.

\begin{figure}[h!]
    \centerline{
    \includegraphics{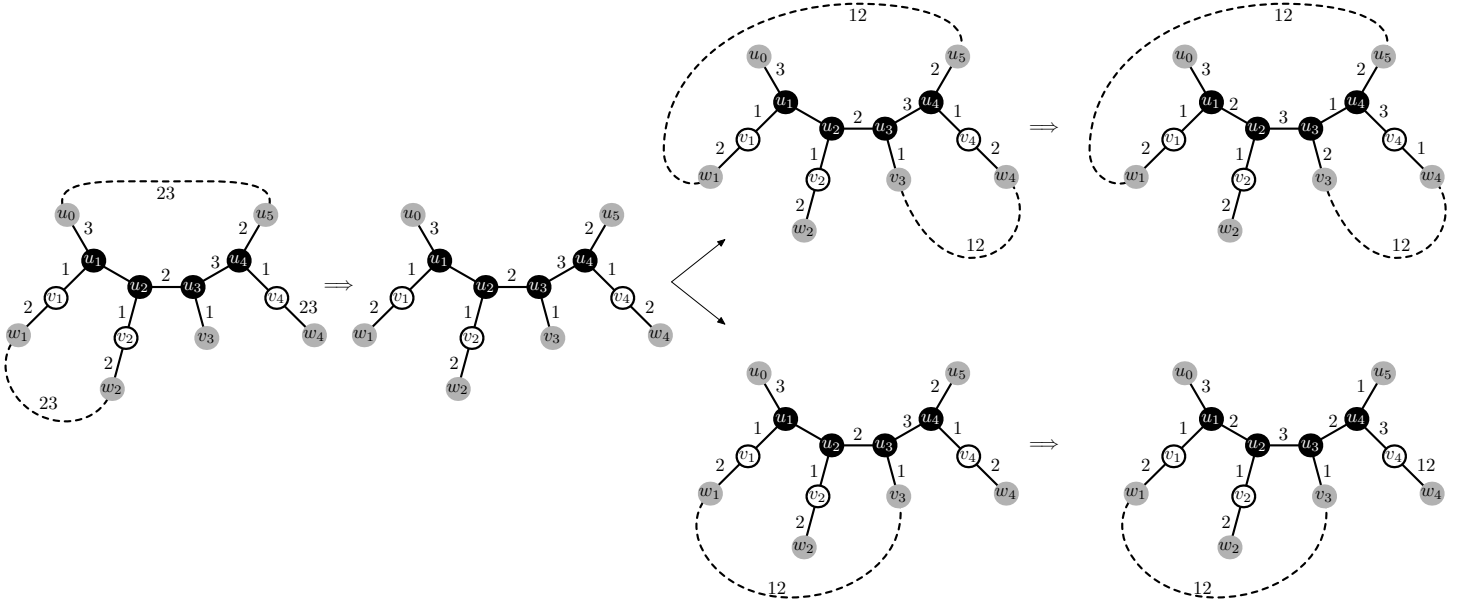}
    }
    \caption{A coloring of $P_{3^2_233_2}$ where the (2,3)-Kempe chain joining $u_1$ and $u_2$ leaves the pattern via the edge $u_4u_5$, is reducible.}
    \label{fig:three0}
\end{figure}

The (1,2)-Kempe chain starting with $u_2u_3$ at $u_2$ still must end with $w_1v_1$ at $v_1$. This chain can contain the path $w_4v_4u_4u_5$ or not. If the chain contains the path, then we can switch only the subchain starting with $u_3v_3$ at $u_3$ and ending with $w_4v_4$ at $v_4$, and then complete the coloring by setting $\varphi(u_1u_2)=2$, $\varphi(u_3u_4)=1$ and $\varphi(v_4u_4)=\varphi(u_2u_3)=3$. Otherwise, we set $\varphi(u_1u_2)=\varphi(u_3u_4)=2$ and $\varphi(u_2u_3)=\varphi(v_4u_4)=3$, after switching the (1,2)-Kempe chain starting with $u_4u_5$ at $u_4$. This chain can eventually reach back the pattern by $w_4v_4$, which does not change the solution. See Figure \ref{fig:three0} for illustration.

Suppose now that the (2,3)-Kempe chain starting with $u_2u_3$ at $u_2$ leaves the pattern by $v_4w_4$. Since the (1,2)-Kempe chain starting with $u_2u_3$ at $u_2$ has to end with $w_1v_1$ at $v_1$, we are free to switch the (1,2)-Kempe chain starting at $v_4$ with $v_4u_4u_5$, and so this case can be reduced to the previous one.

From this point on we may assume that $\varphi(u_3v_3)=3$. To complete the proof, it suffices to reduce this case to the one where we will have $\varphi(u_3v_3)=1$.

Suppose first that the (1,2)-Kempe chain starting with $u_2u_3$ at $u_2$ leaves the pattern by $u_4u_5$. This chain has to end with $w_1v_1$ at $v_1$. There is another (1,2)-Kempe chain starting with $v_4w_4$ at $v_4$; after an eventual switch of the latter we obtain a coloring with $\varphi(v_4w_4)=2$. 

The (2,3)-Kempe chain starting with $u_2u_3$ and $u_3v_3$ at $u_2$ still must end with $u_0u_1$ at $u_1$. This chain can contain the path $w_4v_4u_4u_5$ or not. In both cases we set $\varphi(u_1u_2)=2$, $\varphi(u_2v_2)=3$, and $\varphi(u_2u_3)=\varphi(u_4v_4)=1$.
If the chain contains the path, then we can only switch the subchain starting with $u_3v_3$ at $u_3$ and ending with $w_4v_4$ at $v_4$, and then easily complete the coloring by setting $\varphi(u_3u_4)=3$. Otherwise, we set $\varphi(u_3u_4)=2$ and we switch the (2,3)-Kempe chain starting with $u_4u_5$ at $u_4$. This chain can only eventually reach back the pattern by $w_4v_4$, which does not change the solution.

Suppose now that the (1,2)-Kempe chain starting with $u_2u_3$ at $u_2$ leaves the pattern by $v_4w_4$. Since the (2,3)-Kempe chain starting with $u_2u_3$ at $u_2$ has to end with $u_0u_1$ at $u_1$, we are free to switch the (2,3)-Kempe chain starting at $v_4$ with $v_4u_4u_5$, and so this case can be reduced to the previous one.

\subsection*{The pattern $P_7$}

Let $f=u_1u_2u_3u_4u_5u_6u_7$ be a face of size 7, such that $\degree(u_1)=2$, and let $v_i$ be the third neighbor of $u_i$ for $i=2,\dots 7$. Let $v_4$ and $v_5$ be of degree 2; let $w_i$ be the other neighbor of $v_i$ for $i=4,5$. 

We already know that we may suppose that $\varphi(u_4v_4)=\varphi(u_5v_5)=1$, $\varphi(u_3u_4)=2$, $\varphi(u_5u_6)=3$, and there is a (2,3)-Kempe chain $C$ starting with $u_5u_6$ at $u_5$ and ending with $u_3u_4$ at $u_4$. Moreover, there is another (2,3)-Kempe chain joining $v_4$ and $v_5$, and the colors of $v_4w_4$ and $v_5w_5$ are the same; me may assume that $\varphi(v_4w_4)=\varphi(v_5w_5)=2$. 

If the chain $C$ does not leave the pattern, then $\varphi(u_iv_i)=1$ for $i=2,3,6,7$. The (1,2)-Kempe chain starting with $u_3v_3$ at $u_3$ must end with $w_5v_5$ at $v_5$, after eventually passing by the path $v_7u_7u_6v_6$. After switching $C$, we get a symmetric situation: the (1,2)-Kempe chain starting with $u_6v_6$ at $u_6$ must end with $w_4v_4$ at $v_4$, after eventually passing by the path $v_2u_2u_3v_3$. However, both are not possible due to planarity.

From now on we know that the chain $C$ leaves the pattern. Let us consider the longer subchain, starting either with $u_4u_3$ or with $u_5u_6$, contained in the pattern. By symmetry, after an eventual switch of (2,3)-Kempe chain joining $v_4$ and $v_5$, we may assume that the one that starts with $u_5u_6$ is the longer of the two.

If the (2,3)-Kempe chain starting with $u_5u_6$ takes four consecutive edges on the face $f$, then we can switch the (1,2)-chain $u_4u_3u_2u_1$ and complete the coloring by setting $\varphi(u_4v_4)=3$ and $\varphi(u_4u_5)=2$. The (2,3)-Kempe chain starting with $u_5u_6$ cannot take exactly three consecutive edges on the face $f$, otherwise it could not reach $u_3u_4$.

Let the (2,3)-Kempe chain starting with $u_5u_6$ take exactly two consecutive edges on the face $f$. Then $\varphi(u_1u_7)=1$ and so the other part of $C$, the (2,3)-Kempe chain starting with $u_4u_3$ cannot take exactly two consecutive edges on the face $f$ too. Therefore, $\varphi(u_3v_3)=3$ and $\varphi(u_2u_3)=1$. We switch $C$ and so we have $\varphi(u_3v_3)=\varphi(u_7v_7)=2$. Besides the (2,3)-Kempe chains joining $u_4$ to $u_5$ and $v_4$ to $v_5$, there is a third, disjoint one, starting with $u_1u_2$ and $u_2v_2$ at $u_1$. After an eventual switch, we get $\varphi(u_1u_2)=3$ and $\varphi(u_2v_2)=2$. In other words, all the half-edges leaving the pattern but $u_6v_6$ are colored 2.

\begin{figure}[h!]
    \centerline{
    \includegraphics{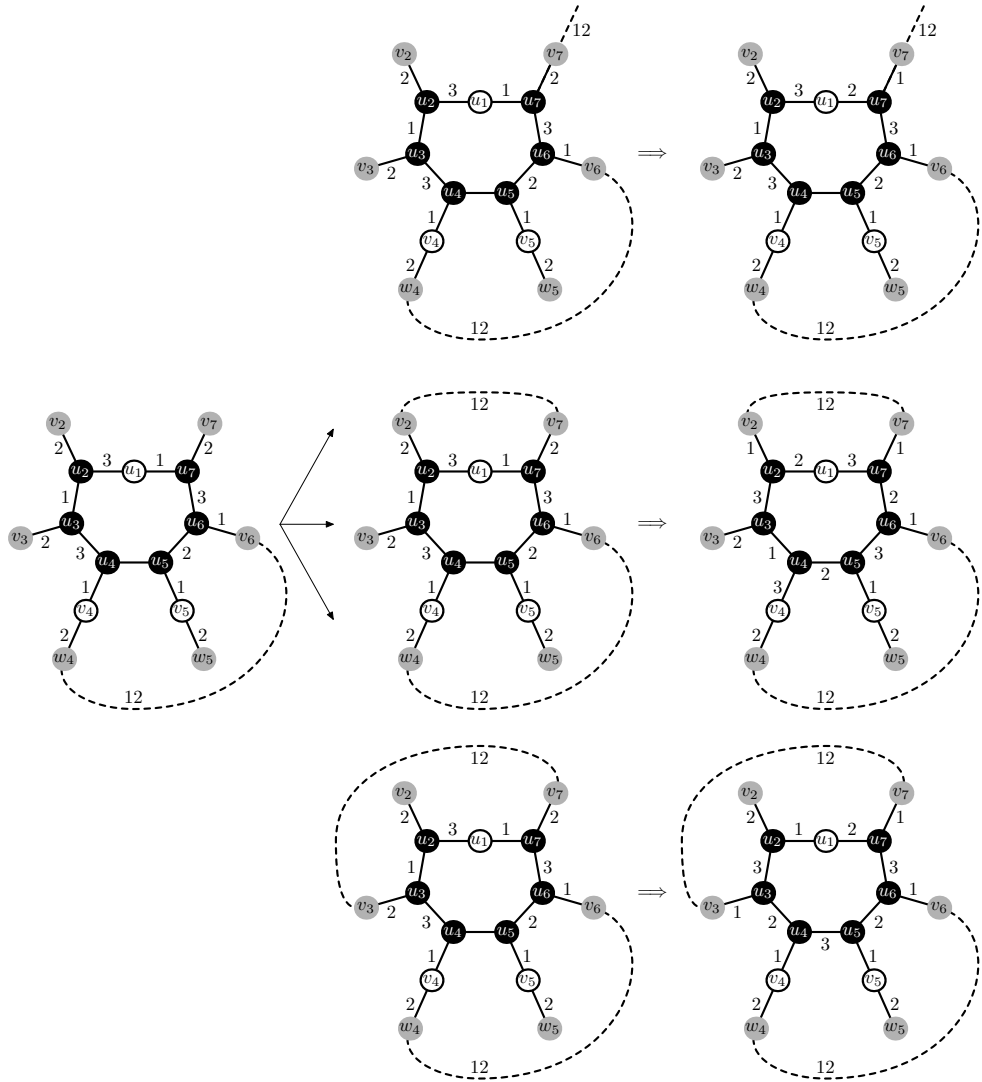}
    }
    \caption{A coloring of $P_7$ where the (1,2)-Kempe chain joining $u_4$ and $u_5$ does not revisit the pattern, is reducible.}
    \label{fig:seven}
\end{figure}

The (1,2)-Kempe chain $C'$ starting with $u_5u_6$ at $u_5$ has to end with $w_4v_4$ at $v_4$. It can eventually contain the path $P'=v_2u_2u_3v_3$. If $C'$ contains $P'$, then we can switch freely the (1,2)-Kempe chain starting at $u_1$ with $u_1u_7$ and $u_7v_7$, reducing to the case where a (2,3)-Kempe chain from $u_5$ takes more that two consecutive edges on the face $f$. If $C'$ does not contain $P'$, then again we consider the (1,2)-Kempe chain $C''$ starting at $u_1$ with $u_1u_7$ and $u_7v_7$. If it does not come back to the pattern, then again we can reduce to the case where $C$ takes four consecutive edges on $f$. The chain $C''$ can come back either by $v_2u_2$ or by $v_3u_3$. In both cases, after a switch of $C''$ we can complete the coloring, see Figure \ref{fig:seven}.

It remains only to consider the case where both ends of $C$, the (2,3)-Kempe chains starting at $u_4$ and at $v_4$, only take one edge on the face $f$. Then $\varphi(u_2u_3)=\varphi(u_6u_7)=1$ and $\{\varphi(u_2v_2),\varphi(u_7v_7)\})=\{2,3\}$. 
If $C$ does not contain the path $v_2u_2u_1u_7v_7$, then after an eventual switch of the Kempe chain containing this path we can get $\varphi(u_2v_2)=\varphi(u_3v_3)$ (and so also $\varphi(u_6v_6)=\varphi(u_7v_7)$), a coloring that is extendible, just like in the case where $C$ had four consecutive edges on $f$.
Hence, we may assume that $C$ contains $v_2u_2u_1u_7v_7$, moreover, $\varphi(u_2v_2)\ne \varphi(u_3v_3)$ and $\varphi(u_6v_6)\ne \varphi(u_7v_7)$. Now the (1,2)-Kempe chain starting with $u_4u_3u_2v_2$ has to end with $w_5v_5u_5$, and so we are free to switch the Kempe chain starting with $u_1u_7u_6v_6$, to reduce this case to the one where at least one branch of $C$ takes at least two consecutive edges on $f$.

\end{document}